\documentclass{gtart_h}

\def\ifplaintex{\expandafter\ifx\csname documentclass\endcsname\relax}

\def\gtp{{\mathsurround=0pt\it $\cal G\mskip-2mu$eometry \&\ 
$\cal T\!\!$opology $\cal P\!$ublications}}  

\def\recd{{\small Received:\qua\receiveddate\ifx\reviseddate\relax
\else\qquad Revised:\qua\reviseddate\fi\par}} 

\def\lognumber#1{\def\thelognumber{#1}}
\def\volumenumber#1{\def\thevolumenumber{#1}}
\def\volumeyear#1{\def\thevolumeyear{#1}}
\def\papernumber#1{\def\thepapernumber{#1}}
\def\pagenumbers#1#2{\def\startpage{#1}\def\finishpage{#2}}
\def\published#1{\def\publishdate{#1}}

\def\received#1{\def\receiveddate{#1}}
\def\revised#1{\def\reviseddate{#1}}
\def\accepted#1{\def\accepteddate{#1}}

\long\def\asciiabstract#1{\long\def\theasciiabstract{#1}}

\let\\\par\let\thelognumber\relax\let\thevolumenumber\relax
\let\thepapernumber\relax\let\thevolumeyear\relax\let\startpage\relax
\let\finishpage\relax\let\publishdate\relax\let\receiveddate\relax
\let\reviseddate\relax\let\accepteddate\relax\let\theasciititle\relax
\let\theasciiauthors\relax
\let\theasciiabstract\relax

\let\theasciiemail\relax

\ifplaintex
\font\logobig=cmssbx10 scaled 3836
\font\logomed=cmssbx10 scaled 2557
\else
\font\logobig=cmssbx10 scaled 4200
\font\logomed=cmssbx10 scaled 2800
\fi

\long\def\makeagttitle{   
\count0=\startpage
\agt\hfill      
\hbox to 45truept{\vbox to 0pt{\vglue -13truept{\logomed A\kern -.37em{\logobig 
T}\kern -.38em G}\vss}\hss}
\break
{\small Volume \thevolumenumber\ (\thevolumeyear)
\startpage--\finishpage\nl
Published: \publishdate}

\vglue .25truein

{\parskip=0pt\leftskip 0pt plus
1fil\def\\{\par\smallskip}{\Large\bf\thetitle}\par\medskip} \vglue
0.05truein

{\parskip=0pt\leftskip 0pt plus 1fil\def\\{\par}{\sc\theauthors}
\par\medskip}%
 
\vglue 0.03truein 

{\small\leftskip 25truept\rightskip 25truept{\bf Abstract}\stdspace\theabstract

{\bf AMS Classification}\stdspace\theprimaryclass
\ifx\thesecondaryclass\relax\else; \thesecondaryclass\fi\par
{\bf Keywords}\stdspace \thekeywords\par}\vglue 7truept

}   

\ifplaintex
\font\phead=cmsl9 scaled 950
\font\pnum=cmbx10 scaled 913
\font\pfoot=cmsl9 scaled 950
\headline{\vbox to 0pt{\vskip -4.5mm\line{\small\phead\ifnum
\count0=\startpage ISSN 1472-2739 (on-line) 1472-2747 (printed)
\hfill {\pnum\folio}\else\ifodd\count0\def\\{ }%
\ifx\theshorttitle\relax\thetitle\else\theshorttitle\fi\hfill{\pnum\folio}
\else\def\\{ and }{\pnum\folio}\hfill\ifx\theshortauthors\relax\theauthors
\else\theshortauthors\fi\fi\fi}\vss}}
\footline{\vbox to 0pt{\vglue 0mm\line{\small\pfoot\ifnum\count0=\startpage
\copyright\ \gtp\hfill\else
\agt, Volume \thevolumenumber\ (\thevolumeyear)\hfill\fi}\vss}}
\else
\font=cmsl9 scaled 1050
\font\lnum=cmbx10 
\font=cmsl9 scaled 1050
\makeatletter
\def\@oddhead{{\small\lhead\ifnum\count0=\startpage ISSN 1472-2739 
(on-line) 1472-2747 (printed)\hfill {\lnum\number\count0}\else\ifodd\count0
\def\\{ }\ifx\theshorttitle\relax \thetitle \else\theshorttitle\fi\hfill
{\lnum\number\count0}\else\def\\{ and }{\lnum\number\count0}
\hfill\ifx\theshortauthors\relax 
\theauthors\else\theshortauthors\fi\fi\fi}}\def\@evenhead{\@oddhead}
\def\@oddfoot{\small\lfoot\ifnum\count0=\startpage\copyright\ \gtp\hfill\else
\agt, Volume \thevolumenumber\ (\thevolumeyear)\hfill\fi}
\def\@evenfoot{\@oddfoot}
\makeatother
\fi
\let\maketitlepage\makeagttitle
\let\makeshorttitle\maketitlepage
\let\maketitle\maketitlepage

\newwrite\gtoutfile
\long\gdef\makeheadfile{  
{\def\\{, }\def\s{ }
\immediate\openout\gtoutfile head.xxx
\immediate\write\gtoutfile{Proxy-for: \ifx\theasciiauthors\relax
\theauthors\else\theasciiauthors\fi\s<\ifx\theasciiemail\relax\theemail\else\theasciiemail\fi>}
\immediate\write\gtoutfile{\noexpand\\}
\immediate\write\gtoutfile{Authors: \ifx\theasciiauthors\relax
\theauthors\else\theasciiauthors\fi}
{\def\\{ }\immediate\write\gtoutfile{Title: \ifx\theasciititle\relax
\thetitle\else\theasciititle\fi}}
\immediate\write\gtoutfile{Subj-class: GT or SG, GR etc}
\immediate\write\gtoutfile{MSC-class: \theprimaryclass\ifx\thesecondaryclass\relax\else, \thesecondaryclass\fi}
\immediate\write\gtoutfile{Journal-ref: Algebr. Geom. Topol. \thevolumenumber\s
(\thevolumeyear) \startpage-\finishpage}
\immediate\write\gtoutfile{Comments: Published by Algebraic and
Geometric Topology at}
\immediate\write\gtoutfile{\s\s\s  http://www.maths.warwick.ac.uk/agt/AGTVol\thevolumenumber/agt-\thevolumenumber-\thepapernumber.abs.html}
\immediate\write\gtoutfile{\noexpand\\}
\immediate\write\gtoutfile{}
\ifx\theasciiabstract\relax
\immediate\write\gtoutfile{\theabstract}\else
\immediate\write\gtoutfile{\theasciiabstract}\fi
\immediate\write\gtoutfile{}
\immediate\write\gtoutfile{\noexpand\\}
\immediate\write\gtoutfile{}
\immediate\closeout\gtoutfile}}  

\def\maketitlepage{\makeagttitle\makeheadfile}
\let\makeshorttitle\maketitlepage
\let\maketitle\maketitlepage

\lognumber{12}
\volumenumber{5}
\volumeyear{2005}
\papernumber{12}
\pagenumbers{219}{235}
\received{17 March 2004} 
\revised{22 March 2005}
\accepted{24 March 2005}
\published{6 April 2005}

\usepackage{amsmath,amssymb,amscd}

\newtheorem{theorem}{Theorem}[section]

 \newtheorem{conj}[theorem]{Conjecture}
\newtheorem{proposition}[theorem]{Proposition}

\newcommand{\edim}{{\rm e-dim}}
\newcommand{\z}{{\Bbb Z}}
\newcommand{\q}{{\Bbb Q}}

\newcommand{\N}{{\Bbb N}}
\newcommand{\s}{{\Bbb S}}
\newcommand {\h}{{\check H}}

\newcommand{\dirlim}{\raisebox{-1ex}{$\stackrel{\hbox{lim}}{\rightarrow}$}}
\newcommand{\invlim}{\raisebox{-1ex}{$\stackrel{\hbox{lim}}{\leftarrow}$}}

\newcommand{\lo}{\longrightarrow}
\newcommand{\sm}{\setminus}
\newcommand{\tor}{{\rm Tor}}

\newenvironment{Relax}{\relax}{\relax}
\begin{document}

\title{Rational acyclic resolutions}

\author{Michael  Levin}
\address{Department of Mathematics, Ben Gurion University of the 
Negev\\P.O.B. 653, Be'er Sheva 84105, ISRAEL}
  
\email{mlevine@math.bgu.ac.il}

\begin{abstract}      
Let $X$ be a compactum such that $\dim_{\q} X \leq n$, $n\geq 2$.  We
prove that there is a ${\q}$-acyclic resolution $r: Z \lo X$ from a
compactum $Z$ of $\dim \leq n$.  This allows us to give a complete
description of all the cases when for a compactum $X$ and an abelian
group $G$ such that $\dim_G X \leq n$, $n \geq 2$ there is a
$G$-acyclic resolution $r : Z \lo X$ from a compactum $Z$ of $\dim
\leq n$.
\end{abstract}
\asciiabstract{%
Let X be a compactum such that dim_Q X <= n, n>1.  We prove that there
is a Q-acyclic resolution r: Z-->X from a compactum Z of dim <= n.
This allows us to give a complete description of all the cases when
for a compactum X and an abelian group G such that dim_G X <= n, n>1
there is a G-acyclic resolution r: Z-->X from a compactum Z of dim <=
n.}

\begin{Relax}\end{Relax}

\primaryclass{55M10, 54F45}
\keywords{Cohomological dimension, acyclic resolution}
\maketitle

\begin{section}{Introduction}

The cohomological dimension $\dim_G X$ of a separable metric space $X$
with respect to an abelian group $G$ is the least number $n$ such that
$\h^{n+1} (X,A;G)=0$ for every closed subset $A$ of $X$.  It was known
long ago that $\dim X =\dim_\z X$ if $X$ is finite dimensional.
Solving an outstanding problem in cohomological dimension theory
Dranishnikov constructed in 1987 an infinite dimensional
compactum (= compact metric space) with $\dim_\z =3$. A few years
earlier a deep relation between $\dim_\z$ and $\dim$ was established
by the Edwards cell-like resolution theorem \cite{ed1,w1} saying that
a compactum of $\dim_\z \leq n$ can be obtained as the image of a
cell-like map defined on a compactum of $\dim \leq n$.  A compactum
$X$ is cell-like if any map $f : X \lo K$ from $X$ to a CW-complex $K$
is null homotopic.  A map is cell-like if its fibers are
cell-like. The reduced $\rm {\check C}$ech cohomology groups of a
cell-like compactum are trivial with respect to any group $G$.

Acyclic resolutions are a natural generalization of the Edwards
cell-like resolution.  A space is $G$-acyclic if its reduced $\rm
{\check C}$ech cohomology groups modulo $G$ are trivial, a map is
$G$-acyclic if every fiber is $G$-acyclic.  As a consequence of the
Vietoris-Begle theorem a surjective $G$-acyclic map of compacta cannot
raise the cohomological dimension $\dim_G$.

Let $\cal P$ denote the set of prime numbers and let ${\cal L} \subset
\cal P$.  Recall that
$${\z_{({\cal L})}} = \{ m/n : m \in \z, n\in \N \text{ and }n\text{ is not
  divisible by the elements of }{\cal L} \}.$$
Thus $ \z_{({\emptyset)}} =\q$ and $ \z_{({\cal P)}} =\z$.
    Dranishnikov \cite{d0.5, d1} proved the following:

\begin{theorem}{\rm \cite{d0.5, d1}}\label{t1}\qua
 Let ${\cal L} \subset \cal P$ and let $X$ be a compactum with
 $\dim_{\z_{({\cal L})}} X \leq n$, $ n \geq 2$.  Then there are a
 compactum $Z $ with $\dim Z \leq n+1$ and $\dim_{\z_{({\cal L})}} Z
 \leq n$ and a ${\z_{({\cal L})}}$-acyclic map $r : Z \lo X$ from $Z$
 onto $X$.
\end{theorem}

In \cite{l1} Theorem \ref{t1} was generalized for arbitrary groups.
      
\begin{theorem}{\rm \cite{l1}}\label{t1.5}\qua
 Let $G$ be an abelian group and let $X$ be a compactum with $\dim_ G
 X \leq n$, $ n \geq 2$.  Then there are a compactum $Z $ with $\dim Z
 \leq n+1$ and $\dim_{G} Z \leq n$ and a $G$-acyclic map $r : Z \lo X$
 from $Z$ onto $X$.
\end{theorem}

It is known that in general in Theorem \ref{t1.5} the dimension of $Z$
cannot be reduced to $n$ \cite{ky2}.  However, Dranishnikov \cite{d1}
showed that it can be done for $G=\z_p$.

 \begin{theorem} {\rm (\cite{d1}, cf.\ \cite{l4})}
\label{t3}
Let ${\cal L} \subset \cal P $ and let $X$ be a compactum with $\dim_
 {\z_p} X\leq n$ for every $p \in \cal L$.  Then there are a compactum
 $Z $ with $\dim Z \leq n$ and a map $r : Z \lo X$ from $Z$ onto $X$
 such that $r$ is ${\z_p}$-acyclic for every $p \in \cal L$.
 \end{theorem}

One of the key problems in the area of acyclic resolutions is whether
 the dimension of $Z$ in Theorem \ref{t1} can be reduced to $n$.  The
 main purpose of this paper is to answer this problem affirmatively by
 proving:

\begin{theorem}
\label{t2}
  Let ${\cal L} \subset \cal P $ and let $X$ be a compactum with
 $\dim_ {\z_{({\cal L})}} X\leq n$, $ n \geq 2$.  Then there are a
 compactum $Z $ with $\dim Z \leq n$ and a ${\z_{({\cal L})}}$-acyclic
 map $r : Z \lo X$ from $Z$ onto $X$.
\end{theorem}

 Theorems \ref{t3} and \ref{t2} allow us to give a complete
 description of all the cases when for an abelian group $G$ and a
 compactum $X$ such that $\dim_G X \leq n$, $n \geq 2$ there is a
 $G$-acyclic resolution $r : Z \lo X$ from a compactum $Z$ of $\dim
 \leq n$.

Recall that the Bockstein basis of abelian groups is the following
 collection of groups $\sigma = \{\q, \z_p, \z_{p^\infty}, \z_{(p)} :
 p \in {\cal P} \}$.  The Bockstein basis $\sigma(G)$ of an abelian
 group $G$ is a subcollection of $\sigma$ defined as follows: 

 {\leftskip 25pt
   $\z_{(p)} \in \sigma(G)$ if $G/\tor G$ is not divisible by $p$;

    $\z_p \in \sigma(G)$ if $\tor_p G $ is not divisible by $p$;

    $\z_{p^\infty} \in \sigma(G)$ if $\tor_p G \neq 0$ and $\tor_p G $
    is divisible by $p$;

    $ \q \in \sigma(G)$ if $G/\tor G\neq 0$ and $G/\tor G$ is
divisible by every $p \in \cal P$.

}  

Let $X$ be a compactum. The Bockstein theorem says that $$\dim_G X =
\sup\{\dim_E X : E \in \sigma(G) \}.$$  It is known \cite{ky2} that $X$
is $G$-acyclic if and only if $X$ is $E$-acyclic for every $E \in
\sigma(G)$.

 Following \cite{l5} define the closure $\overline{\sigma({ G})}$ of
the Bockstein basis $\sigma({ G}) $ of $G$ as a collection of abelian
groups containing ${\sigma({ G})}$ and possibly some additional groups
determined by:

 {\leftskip 25pt
$\z_p \in \overline{\sigma({ G})}$ if $\z_{p^\infty} \in
\overline{\sigma({ G})}$;

$\z_{p^\infty} \in \overline{\sigma({ G})}$ if $\z_{p} \in
\overline{\sigma({ G})}$;

$\q \in \overline{\sigma({G})}$ if $\z_{(p)} \in \overline{\sigma({
G})}$;

$\z_{(p)} \in \overline{\sigma({G})}$ if $\q $ and $\z_{p^\infty} \in
\overline{\sigma({G})}$.\par}

\begin{theorem}
\label{t4} Let $G$ be an abelian group, let $X$ be a compactum
and let $n\geq 2$.  Then there exist a compactum $Z$ of $\dim \leq n$
 and a $G$-acyclic map $r : Z \lo X$ from $Z$ onto $X$ if and only if
 $\dim_E X \leq n$ for every $E \in \overline{\sigma({ G})}$.
\end{theorem}

\proof
Assume that there exist a compactum $Z$ of $\dim \leq
   n$ and a $G$-acyclic map $r : Z \lo X$ from $Z$ onto $X$.  Then $r$
   is $E$-acyclic for every $E \in \sigma({ G})$.  Consider the exact
   sequences:

 {\leftskip 25pt
    $0 \lo \z_p \lo \z_{p^\infty} \lo\z_{p^\infty}\lo 0$;

    $0 \lo \z_{p^k} \lo \z_{p^{k+1}} \lo \z_p \lo 0$;

    $0 \lo \z_{(p)} \lo \z_{(p)} \lo \z_p \lo 0$;

    $0 \lo \z_{(p)} \lo \q \lo \z_{p^\infty} \lo 0$.\par}  

Then it follows from the Bockstein sequences generated by the above
 exact sequences that in the class of compacta:

 {\leftskip 25pt
the $ \z_{p^\infty}$-acyclicity implies the $\z_p $-acyclicity;

the $\z_p $-acyclicity implies the $ \z_{p^k}$-acyclicity for every
 $k$ and since $\z_{p^\infty} =\dirlim \z_{p^k}$ the $\z_p
 $-acyclicity implies the $\z_{p^\infty} $-acyclicity;

the $ \z_{(p)}$-acyclicity implies the $\z_p $-acyclicity;

the $ \z_{(p)}$-acyclicity implies the $\q$-acyclicity;

the $\q$-acyclicity together with the $ \z_{p^\infty}$-acyclicity
imply the $ \z_{(p)}$-acyclicity.\par}  

From these properties we obtain that $r$ is $E$-acyclic for every $E
\in \overline{\sigma({G})} $ and therefore by the Vietoris-Begle
theorem $\dim_E X \leq n$ for every $E \in \overline{\sigma({ G})} $.

Now assume that $\dim_E X \leq n $ for every $E \in
\overline{\sigma({G})} $.

Consider first the case when $G$ is torsion.  Let ${\cal L}=\{
 p\in{\cal P} : \tor_p G \neq 0 \}$.  Then $\dim_{\z_p} X \leq n$ for
 every $p \in \cal L$ and the existence of the required resolution
 follows from Theorem \ref{t3} and the above properties of acyclicity.

Assume that $G$ is not torsion.  Then $\q \in \overline{\sigma({G})}
$.  Define ${\cal L} =\{ p \in {\cal P}: \z_{(p)} \in
\overline{\sigma({G}) }\}$.  By the Bockstein theorem
$\dim_{\z_{({\cal L})} }X \leq n$ and therefore by Theorem \ref{t2}
there exists a $\z_{({\cal L}) }$-acyclic resolution $r : Z \lo X$
from a compactum $Z$ of $\dim \leq n$ onto $X$.  Clearly $r$ is
$\q$-acyclic if ${\cal L} =\emptyset$.  If ${\cal L} \neq \emptyset$
then $r$ is $ \z_{(p)}$-acyclic for every $\z_{(p)} \in
\overline{\sigma({ G})} $ and by the above properties $r$ is
$\q$-acyclic as well.

 Let $ \z_{(p)} \in \sigma( G)$. Then $ \z_{(p)} \in
 \overline{\sigma({G})} $ and therefore $r$ is $ \z_{(p)}$-acyclic.

  Let $\z_p $ or $\z_{p^\infty} \in \sigma(G)$.  Then $ \z_{p^\infty}
  \in \overline{\sigma({ G})} $ and since $\q \in \overline{\sigma({
  G})}$ we have that $\z_{(p)} \in \overline{\sigma({ G}) }$. Hence
  $r$ is $\z_{(p)}$-acyclic and by the above properties $r$ is also
  $\z_p $-acyclic and $\z_{p^\infty}$-acyclic.

Thus we have shown that $r$ is $E$-acyclic for every $ E \in
\sigma(G)$ and hence $r$ is $G$-acyclic.  The theorem is proved.
\endproof

Theorem \ref{t1.5} was generalized in \cite{l5} as
follows.

 \begin{theorem}{\rm \cite{l5}}\label{t5}\qua
Let $X$ be a compactum.  Then for every integer $n \geq 2$ there are a
  compactum $Z$ of $\dim \leq n+1$ and a surjective map $r: Z \lo X$
  having the property that for every abelian group $G$ and every
  integer $k \geq 2$ such that $\dim_G X \leq k \leq n$ we have that
  $\dim_G Z \leq k$ and $r$ is $G$-acyclic.  \end{theorem} 

Theorem \ref{t5} resulted in the following conjecture posed in
  \cite{l5}.
  
\begin{conj}{\rm\cite{l5}}\label{con}\qua
Let ${\cal G}$ be a collection of abelian groups, let $X$ be a
  compactum $X$ and let $n \geq 2$.  Then there exist a compactum $Z$
  of $\dim \leq n$ and a surjective map $r : Z \lo X$ such that
  $\dim_G Z \leq \max\{\dim_G X , 2\}$ and $r$ is $G$-acyclic for
  every $G \in \cal G$ if and only if $\dim_E X \leq n$ for every $E
  \in \overline{\sigma(\cal G)}=\cup\{ \overline {\sigma(G)} :G \in
  {\cal G}\}$.
\end{conj}

It is easy to derive from the proof of Theorem \ref{t4} that the
  condition $\dim_E X \leq n$ for every $E \in \overline{\sigma(\cal
  G)}$ is necessary for the existence of a resolution required in the
  conjecture.  Theorem \ref{t4} seems to be an important step toward
  proving that this condition is also sufficient.  

 Note that neither the restriction $k \geq 2$ in Theorem \ref{t5} nor
 the restriction $\dim_G Z \leq \max\{\dim_G X , 2\}$ in Conjecture
 \ref{con} can be replaced by $k \geq 1$ and $\dim_G Z \leq \dim_G X$
 respectively.

 Indeed, let a compactum $Y$ be $\q$-acyclic and $\dim_\q Y \leq 1$.
 Then for every closed subset $A$ of $Y$ the exact sequence of the
 pair $(Y,A)$ implies that $\h^1(A; \q)$ $=0$ and by Universal
 Coefficient Theorem $\h^1(A) =0$.  Hence $\h^1(Y) =0$ and $\dim Y
 \leq 1$ since every map from $A$ to $\s^1=K(\z,1)$ is null-homotopic
 and therefore extends over $Y$.  Then $Y$ is $\z$-acyclic.  Thus we
 obtain that a $\q$-acyclic map from a compactum of $\dim_\q \leq 1$
 is $\z$-acyclic and hence there is no $\q$-acyclic resolution $r: Z
 \lo X$ for a compactum $X$ of $\dim_\z =\infty$ and $\dim_\q =1$ from
 a finite dimensional compactum $Z$ of $\dim_\q \leq 1$.

 The last observation also shows that the restriction $n \geq 2$ in
 Theorem \ref{t2} cannot be omitted.

 \end{section} \begin{section}{Preliminaries} The extensional
 dimension of a compactum $X$ is said not to exceed a CW-complex $K$,
 written $\edim X \leq K$, if for every closed subset $A$ of $X$ and
 every map $f : A \lo K$ there is an extension of $f$ over $X$.  It is
 well-known that $\dim X \leq n$ is equivalent to $\edim X \leq \s^n$
 and $\dim_G X \leq n$ is equivalent to $\edim X \leq K(G,n)$, where
 $K(G,n)$ is an Eilenberg-Mac Lane complex of type $(G,n)$.  Also note
 that if $\dim X \leq n+1$, $n \geq 2$ then $\dim_G X \leq n$ is
 equivalent to $\edim X \leq M(G,n)$, where $M(G,n)$ is a
 simply-connected Moore complex of type $(G,n)$ (the condition $\dim X
 \leq n+1$ can be relaxed to $\dim X < \infty $ \cite{dr3}).

     A map between CW-complexes is said to be combinatorial if the
   preimage of every subcomplex of the range is a subcomplex of the
   domain.

  Let $M$ be a CW-complex.  We say that a triangulation of $M$ is
  compatible with the CW-structure if every cell is a simplicial
  subcomplex.  

  Let ${\cal L} \subset \cal P$.  We will call a simplicial complex
  $M$ an $n$-dimensional $\cal L$-sphere if the following conditions
  are satisfied:

 {\leftskip 25pt  
  $M$ is a finite $n$-dimensional simplicial complex;

  $\tilde{H_i}( M; \z_{({\cal L})})=0$ if $i \neq n$ and $\tilde{H_n}(
    M; \z_{({\cal L})})=\z_{({\cal L})}$;

    for every $n$-simplex $\Delta$ of $M$, the subcomplex $B=(M \sm
    \Delta)\cup \partial \Delta$ is $\z_{({\cal L})}$-acyclic. (Note
    that by Universal Coefficient Theorem the $\z_{({\cal
    L})}$-acyclicity of $B$ is equivalent to $\tilde{H}_*(B;
    \z_{({\cal L})})=0$.)\par}  

 It is easy to verify that if $M$ is an $n$-dimensional $\cal
  L$-sphere with respect to a given triangulation then $M$ is also an
  $n$-dimensional $\cal L$-sphere with respect to any barycentric
  subdivision of the triangulation.  

 By an $n$-dimensional $\cal L$-ball $B$ we mean an $n$-dimensional
  $\cal L$-sphere without the interior of one of its $n$-simplexes.
  The boundary of the simplex whose interior is removed is said to be
  the boundary of the $\cal L$-ball $B$ and it is denoted by $\partial
  B$.  

 Saying that an $n$-dimensional simplicial complex $N$ is obtained
  from an $n$-dimensional simplicial complex $M$ by replacing an
  $n$-simplex $\Delta $ of $M$ by an $n$-dimensional $\cal L$-ball $B$
  we mean that the interior of $\Delta$ is removed from $M$ and $B$ is
  attached to the boundary of $\Delta$ by a simplicial homeomorphism
  of $ \partial \Delta$ and $\partial B$.  We regard $N$ as a
  simplicial complex with respect to the natural triangulation induced
  by the triangulations of $M$ and $B$.

It is clear that replacing an $n$-simplex of an $n$-dimensional $\cal
L$-sphere by an $n$-dimensional $\cal L$-ball we obtain again an
$n$-dimensional $\cal L$-sphere and replacing an $n$-simplex of an
$n$-dimensional $\cal L$-ball by an $n$-dimensional $\cal L$-ball we
obtain again an $n$-dimensional $\cal L$-ball.

  \begin{proposition} \label{p1} Let ${\cal L} \subset \cal P$, let
  $X$ be an compactum with $\dim_{\z_{({\cal L})}} X \leq n $ and
  $\dim X \leq n+1$, $n\geq 2$ and let $f: X \lo K$ be a map from $X$
  to an $(n+1)$-dimensional finite simplicial complex $K$.  Then for
  every $\epsilon >0$ there exist an $(n+1)$-dimensional finite
  CW-complex $M$ and maps $\phi: X \lo M$, $\psi : M\lo K$ such that:

{\rm(i)}\qua $f$ and $\psi \circ \phi$ are $\epsilon$-close, that is ${ \rho}
(f(x), \psi(\phi(x))) < \epsilon $ for every $x \in X$ where $ {\rho}$
is a metric in $K$;

{\rm(ii)}\qua $\psi$ is combinatorial; $\psi $ is 1-to-1 over every simplex of
 $K$ which is not contained in an $(n+1)$-simplex of $K$; the
 $(n+1)$-cells of $M$ are exactly the preimages of the
 $(n+1)$-simplexes of $K$ under $\psi$ such that the interior points
 of the $(n+1)$-cells of $M$ are sent to the interior of the
 corresponding $(n+1)$-simplexes of $K$;

{\rm(iii)}\qua each $(n+1)$-cell $C$ of $M$ has a point $b \in \partial C$ such
that $b$ is sent by $\psi$ to the interior of the corresponding
$(n+1)$-simplex of $K$.  We will call $b$ a free boundary point of
$C$.  It is clear from (ii) that $C$ is the only $(n+1)$-cell of $M$
containing $b$;

{\rm(iv)}\qua $M$ admits a compatible triangulation with respect to which the
     boundary of every $(n+1)$-cell of $M$ is an $n$-dimensional $\cal
     L$-sphere.  Moreover: by (iii) such a triangulation of $M$ can be
     chosen such that for each cell $C$ of $M$ there is an $n$-simplex
     contained in $\partial C$ and consiting only of free boundary
     points of $C$.

  \end{proposition}

\proof
Fix an $n$-simplex $\Delta$ of $K$ contained in at
least one $(n+1)$-simplex and let $\Delta_1, \Delta_2, ..., \Delta_k$,
$k \geq 1$ be the $(n+1)$-simplexes of $L$ containing $\Delta$.  Take
a small closed $n$-dimensional ball $B$ contained in $\Delta$ such
that $B$ does not touch the boundary of $\Delta$ and $B$ is centered
at the barycenter $c$ of $\Delta$.  For every $i=1,2,..,k$ take a
point $p_i $ sufficiently close to $B$ and contained in the interior
of $\Delta_i$.  Denote $P=\{ p_1,..., p_k\} $.  Consider the join
$B*p_i$ as the subset of $\Delta_i$ which is the cone over $ B$ with
the vertex at $p_i$.  Thus we regard $B*P =\cup B*p_i$ as a subset of
$\cup \Delta_i$.  Let the $(n-1)$-dimensional sphere $S $ be the
boundary $\partial B$ of $B$. Denote $F=f^{-1}( B*P)$ and $A= f^{-1}(
S*P)$.

Represent a Moore space of type $( \z_{{(\cal L})}, n-1)$ as an
infinite telescope $T$ of a sequence of $(n-1)$-dimensional spheres
with bonding maps of non-zero degrees which are not divisible by the
elements of $\cal L$.  By a finite subtelescope of $T$ we mean a
subspace of $T$ which is the telescope of finitely many consecutive
spheres in the sequence.  Define the group $G$ as $G=0$ if $k=1$ and
$G=$ the direct sum $\oplus \z_{{(\cal L})}$ of $k-1$ copies of
$\z_{{(\cal L})}$ if $k \geq 2$. Let $D=\{d_1,d_2,...,d_k\}$ be a
discrete space of $k$ singletons. Then $ T * D$ is a simply connected
Moore space of type $(G,n)$.  Let $\tau: S \lo T$ be an embedding of
$S$ as the first sphere of the telescope $T$ and let $\delta : P \lo
D$ be the map sending $p_i$ to $d_i$.

Consider the map $ \alpha=(\tau*\delta)\circ f |_{...}  : A \lo T*D$.
By Bockstein Theorem $\dim_G X \leq n$ and since $\dim X \leq n+1$ we
have that $\edim X \leq T*D$.  Then $ \alpha : A \lo T*D$ extends to
$\beta : F \lo T*D$.  Since $\beta(F)$ is compact take a finite
subtelescope $T_\beta $ of $T$ starting at the first sphere and ending
at a sphere $S_\beta$ such that $\beta(F) $ is contained in $T_\beta *
D$ and let $\gamma : T_\beta \lo S_\beta $ be the natural retraction
sending $ T_\beta$ to the last sphere of $ T_\beta$.

Define a CW-complex $M_\Delta$ as the quotient space of $L =(K \sm
 (B*P))\cup (S*P)$ obtained by identifying the points of $S*P$ with
 $S_\beta * D$ according to the map $(\gamma \circ \tau)* \delta : S*P
 \lo S_\beta * D$.  Denote by $\pi : L \lo M_\Delta$ the projection
 and consider $ S_\beta * D$ as the subset $\pi(S*P)$ of $M_\Delta $.
 We will also consider $L \sm (B*P)$ as a subset of $ M_\Delta$.  Note
 that since the identifications on $ S_\beta * D$ are induced by the
 join $(\gamma \circ \tau)* \delta$ we have that $\pi^{-1}(\pi(S*p_i))
 = \pi^{-1}(S_\beta *d_i)= S*p_i$, $\pi^{-1}(\pi(p_i)) =
 \pi^{-1}(d_i)= p_i$ and $\pi^{-1}(\pi(S)) = \pi^{-1}(S_\beta )= S$.

  It is easy to see that $ (\gamma* id_D)( \alpha (x))=\pi (f(x)) \in
 S_{\beta} * D$ for $x\in A$ where $id_D : D \lo D$ is the identity
 map.  Then the map $\pi\circ f |_{...}  : (X \sm F) \cup A \lo
 M_\Delta$ extends to the map $\phi_\Delta : X \lo M_\Delta$ defined
 by $\phi_\Delta (x)= (\gamma* id_D)( \beta (x))$ for $x \in F$.  

Define a map $\psi_\Delta : M_\Delta \lo K$ such that:

 {\leftskip 25pt 
 every simplex of $K$ not intersecting $B*P$ and considered as a
 subset of $M_\Delta$ is sent by $\psi_\Delta$ to itself by the
 identity map;

 $\psi_\Delta$ sends $\pi (\Delta \cap L)$ onto $\Delta $ such that
  $\pi(S)=S_\beta$ is sent to the barycenter $c$ of $\Delta$;

  $\psi_\Delta$ sends $\pi(\Delta_i \cap L)$ onto $\Delta_i$ such that
  $\pi((\Delta_i \sm \partial \Delta_i)\cap L)$ is sent into $\Delta_i
  \sm \partial \Delta_i$;

   $\psi_\Delta \circ \phi$ is sufficiently close to $f$ provided that
   the ball $B$ and the points $p_ i \in \Delta_i $ were chosen to be
   sufficiently close to the barycenter $c$ of $\Delta$. \par} 

Then for every $\Delta_i$, $\psi_\Delta^{-1}( \Delta_i)=\pi( \Delta_i
   \cap K)$ and since $\Delta_i \cap K = (\Delta_i \sm (B*p_i))\cup
   (S*p_i)$ is homeomorphic to an $(n+1)$-dimensional Euclidean ball
   in which $\pi$ makes identifications only on the boundary we may
   endow $M_\Delta$ with a CW-structure so that $\psi_\Delta$ becomes
   a combinatorial map and the $(n+1)$-cells of $M_\Delta$ are exactly
   the preimages of the $(n+1)$-simplexes of $K$.  Moreover, since the
   identifications on $S*p_i$ are induced by the map $\gamma : S \lo
   S_\beta$ of degree not divisible by the elements of $\cal L$ one
   can define a compatible triangulation of $M_\Delta$ so that the
   boundary of each cell $C_i=\psi_{\Delta}^{-1}(\Delta_i)$ becomes an
   $\cal L$-sphere.  Finally note that $\pi(p_i)$ is a free boundary
   point of the cell $C_i$ and if a triangulation of $M_i$ is chosen
   to be sufficiently small then for every $i$ there will be an
   $n$-simplex contained in $\partial C_i$ consisting of free boundary
   points of $C_i$.

It is easy to see that the procedure of constructing $M_\Delta$ and
$\psi_\Delta : M_\Delta \lo K$ described above can be carried out
independently for all $n$-simplexes $\Delta$ of $L$ contained in some
$(n+1)$-simplexes of $K$.  This way we can construct a CW-complex $M$
and a map $\psi : M \lo K$ satisfying the conclusions of the
proposition.  \endproof

  Let $K'$ be a simplicial complex and let $\kappa : K \lo K'$,
 $\lambda : L \lo L'$, $\alpha : L \lo K$ and $\alpha' : L' \lo K'$ be
 maps.
 \[ \begin{CD} L @>\alpha>> K\\ @V\lambda VV @V\kappa VV \\ L'
    @>\alpha'>> K' \end{CD}
\]
We say that $\kappa, \lambda ,\alpha$ and $\alpha'$ combinatorially
 commute if for every simplex $\Delta$ of $K'$ we have that $(\alpha'
 \circ\lambda )(( \kappa\circ \alpha)^{-1}(\Delta)) \subset \Delta$.
 (The direction in which we want the maps $ \kappa,\lambda, \alpha$
 and $\alpha'$ to combinatorially commute is indicated by the first
 map in the list. Thus saying that $\alpha', \kappa, \lambda $ and
 $\alpha$ combinatorially commute we would mean that $(\kappa\circ
 \alpha)(( \alpha'\circ \lambda)^{-1}(\Delta)) \subset \Delta$ for
 every simplex $\Delta $ of $K'$.)

   Let $K$ and $K'$ be simplicial complexes and let $\kappa : K \lo
    K'$ $\kappa': K \lo K'$.  We say that $\kappa' $ is a simplicial
    approximation of $\kappa$ if $\kappa'$ is a simplicial map and for
    every simplex $\Delta$ of $K'$, $\kappa'(\kappa^{-1}(\Delta))
    \subset \Delta$.  Thus if $\kappa' $ is a simplicial approximation
    of $\kappa$ then for every subcomplex $ N$ of $K'$ and for every
    $A \subset K$ such that $\kappa(A) \subset N$ we also have
    $\kappa'(A) \subset N$.

 For a simplicial complex $K$ and $a \in K$, $st(a)$ (the star of $a$)
  denotes the union of all the simplexes of $K$ containing $a$.

 The following proposition whose proof is left to the reader is a
 collection of simple combinatorial properties of maps.

 \begin{proposition} \label{p2}${}$

 {\rm(i)}\qua Let $\kappa : K \lo K'$ be a combinatorial map of simplicial
 complexes $K$ and $K'$. Then $\kappa$ admits a simplicial
 approximation.

 {\rm(ii)}\qua Let $\kappa : K \lo K'$ be a map of simplicial complexes $K$ and
 $K'$.  Then there is a sufficiently small barycentric subdivision of
 the triangulation of $K$ with respect to which $\kappa$ admits a
 simplicial approximation.

{\rm(iii)}\qua Let $K$ and $K'$ be simplicial complexes, let maps $\kappa : K
 \lo K'$, $\lambda : L \lo L'$, $\alpha : L \lo K$ and $\alpha' : L'
 \lo K'$ combinatorially commute and let $\kappa$ be combinatorial.
 Then
  $$\lambda( \alpha^{-1}(st( a)) ) \subset \alpha'{}^{-1}(st( \kappa
   (a)) )\text{ and }\kappa(st((\alpha (b))) \subset st((\alpha' \circ
   \lambda)(b))$$
for every $a \in K$ and $b\in L$.
 \end{proposition}

  \end{section}

  \begin{section}{Proof of Theorem \ref{t2}} 

By the Vietoris-Begle theorem the composition of $\z_{({\cal
L})}$-acyclic maps of compacta is again a $\z_{({\cal L})}$-acyclic
map. Then by virtue of Theorem \ref{t1} we may assume without loss of
generality that $\dim X = n+1$.

 For every $i=1,2,...$ we will construct spaces $L_i, W_i, S_i , M_i,
 K_i$, maps represented in the diagram
 \[ \begin{CD} L_{i+1} @>\beta_{i+1}>> W_{i+1} @.\subset S_{i+1}
    \subset @. M_{i+1}@>\psi_{i+1}>> K_{i+1}\\ @V\lambda_{i+1}VV
    @V\omega_{i+1}VV @.@V\mu_{i+1}VV @V\kappa_{i+1}VV \\ L_{i}
    @>\beta_{i}>> W_{i}@. \subset S_i \subset @.M_{i}@>\psi_{i}>>
    K_{i} \end{CD}
\]
and some additional maps including $\phi_i: X \lo M_i$ and $f_i : X
\lo K_i$ having the following properties.  

\medskip
{\bf Remark}\qua Since the list of properties (steps of the construction)
is rather long we would like to outline a way to get the main idea of
the construction omitting some details.  The property (9) is the heart
of the construction.  The fastest way to get to (9) is to skip (8) and
to omit all metric estimates in (1), (2), (6).  Also note that the
main purpose of (8), ( 11), (12) is to show that the maps
$\kappa_{i+1}, \alpha_{i+1}, \lambda_{i+1}$ and $ \alpha_i$
combinatorially commute where $\alpha_{i}=\psi_{i}\circ \beta_{i} :
L_i \lo K_i$.

\medskip
{\bf(1)}\qua $K_i$ is a finite simplicial complex of $\dim\leq n+1$.  We fix
  a metric $d_i$ in $K_i$ such that the distance between any two
  points of $K_i$ lying in non-intersecting simplexes of $K_i$ is
  bigger than $1$.  We also fix a metric in $X$.  

\medskip
{\bf(2)}\qua Define $\epsilon_i$ to be so small that for every $j \leq i$ and
every $x,y \in K_i$ with $d_i(x,y) <\epsilon_i$ we have that
$d_j(\kappa_i^j(x),\kappa_i^j(y)) < 1/4^i$ where $\kappa_i^i
=id_{K_i}$ and $\kappa_i^j =\kappa_{j+1}\circ ...\circ\kappa_ i : K_i
\lo K_j$ for $j < i$.  By Proposition \ref{p1} we choose a CW-complex
$M_i$ and maps $\psi_i : M_i \lo K_i$ and $\phi_i : X \lo M_i$ such
that the conclusions of Proposition \ref{p1} are satisfied with $K$,
$M$, $\phi$, $\psi$, $f$, $\rho$ and $\epsilon$ replaced by $K_i$,
$M_i$, $\phi_i$, $\psi_i$, $f_i$, $d_i$ and $\epsilon_i/2$
respectively.  Thus $M_i$ is a finite CW-complex of $\dim \leq n+1$
admitting a compatible triangulation.  We fix a compatible
triangulation of $M_i$ and consider $M_i$ as both a CW-complex and a
simplicial complex.  Note that the CW-structure of $M_i$ is different
from the simplicial structure of $M_i$ and therefore the cells of
$M_i$ and simplexes of $M_i$ are not the same.  The spaces $M_i$ are
the only spaces in which we will consider simultaneously two
structures.  

  \medskip
{\bf(3)}\qua Let $W_i =$ the $n$-skeleton of $M_i$ with respect to the
      CW-structure of $M_i$.  

 Consider $W_i$ as a simplicial
      complex with a triangulation of $W_i$ to be a sufficiently small
      barycentric subdivision of the triangulation induced by the
      triangulation of $M_i$. We refer to this triangulation of $W_i$
      considering simplexes of $W_i$.  Thus the simplexes of $W_i$ and
      the simplexes of $M_i$ contained in $W_i$ are not the same.
      However each simplicial subcomplex of $M_i$ contained in $W_i$
      will be also a subcomplex of $W_i$.  

 \medskip
{\bf(4)}\qua  Let $S_i=$ the $n$-skeleton of $M_i$ with respect to the
  triangulation of $M_i$.  

 Clearly $W_i \subset S_i$.  Consider
  $S_i$ as a simplicial complex with the triangulation induced by the
  triangulation of $M_i$.

 Let $\gamma_i : S_i \lo W_i$ be a natural retraction built by picking
  up in each $(n+1)$-cell of $M_i$ an interior point not belonging to
  $S_i$ and retracting the complement of the point inside the cell to
  the boundary of the cell.  Then by (ii) of Proposition \ref{p1}
  $$\gamma(\psi_i^{-1}(\Delta) \cap S_i) \subset \psi_i^{-1}(\Delta)\text{
   for every simplex }\Delta\text{ of }K_i.$$  
Indeed, it is obvious if
   $\Delta$ is an $(n+1)$-simplex since then $\psi_i^{-1}(\Delta)$ is
   an $(n+1)$-cell of $M_i$.  If $\Delta$ is of $\dim \leq n$ then
   $\psi_i^{-1}(\Delta) \subset W_i$ and we simply have
   $\gamma(\psi_i^{-1}(\Delta)) =\psi_i^{-1}(\Delta)$.

  It is also easy to see that assuming that the triangulation of $M_i$
  is sufficiently small $\gamma_i$ can be constructed such that for
  every $c \in W_i$ there exists a contractible subcomplex $N$ of
  $M_i$ contained in $W_i$ such that
   $$ \gamma_i(st(c, M_i)\cap S_i) \subset N \subset
   \psi_i^{-1}(st(\psi_i(c)))$$
 where $st(c, M_i)$ is the star of
   $c$ with respect to the triangulation of $M_i$.

   As we already noted in (3), $N$ is also a subcomplex of $W_i$.

 \medskip
{\bf(5)}\qua $L_i$ and $\beta_i : L_i \lo W_i$ are constructed as follows.
  $L_i$ is a simplicial complex of $\dim \leq n$ obtained from $ W_i$
  by replacing some $n$-simplexes of $ W_i$ by $n$-dimensional $\cal
  L$-balls.  Then $\beta_i$ is a projection of $L_i$ to $W_i$ such
  that $\beta_i$ takes the $\cal L$-balls to the corresponding
  replaced simplexes such that the non-boundary points of the $\cal L
  $-balls are sent to the interior of the replaced simplexes and
  $\beta_i$ is 1-to-1 over all the simplexes of $W_i$ which are not
  replaced;

Let $i=1$.  Set $\epsilon_1=1$.  Let $K_1$ be any $(n+1)$-dimensional
simplicial complex with any map $f_1 : X \lo K_1$.  Applying
Proposition \ref{p1} define for $i=1$ all needed spaces and maps
satisfying the relevant properties of (2-4).  Define $L_1$ as $W_1$
with no simplexes replaced by $\cal L$-balls.  We proceed from $i$ to
$i+1$ as follows.  

\medskip
{\bf(6)}\qua Based on $\phi_i : X \lo M_i$ find a finite simplicial complex
  $K_{i+1}$ and maps $f_{i+1} : X \lo K_{i+1}$ and $\nu_{i+1} :
  K_{i+1} \lo M_i$ such that:

{\leftskip 25pt
$d_i((\psi_i \circ \nu_{i+1}\circ f_{i+1}) (x), ( \psi_i\circ
\phi_i)(x)) < \epsilon_i/2$ for every $x\in X$ with $\epsilon_i$
defined in (1);

$f_{i+1}^{-1}(\Delta)$ is of diam $\leq 1/(i+1)$ for every simplex
$\Delta$ of $K_{i+1}$;

$\nu_{i+1}$ is a combinatorial map with respect to the triangulation
 of $M_i$.\par}  

 Define $\kappa_{i+1}= \psi_i \circ \nu_{i+1}: K_{i+1}
 \lo K_i$ and suppose that the triangulation of $K_{i+1}$ is so small
 that

{\leftskip 25pt
$\kappa_{i+1}^j( \Delta)$ is of diam$\leq 1/(i+1)$ for every simplex
  $\Delta$ of $K_{i+1}$ and $j \leq i$.\par}  

\medskip
{\bf(7)}\qua Let $M_{i+1},S_{i+1}$ , $W_{i+1}$ and the corresponding maps be
 defined as described in (2-4).  Since $ \nu_{i+1}\circ\psi_{i+1} :
 M_{i+1} \lo M_i$ is combinatorial with respect to the simplicial
 structures of $M_{i+1}$ and $M_i$, by (i) of Proposition \ref{p2}
 there is a simplicial approximation $\mu_{i+1} : M_{i+1} \lo M_i$ of
 $ \nu_{i+1}\circ\psi_{i+1}$.  Then $\mu_{i+1}(S_{i+1}) \subset S_i$.
 By (ii) of Proposition \ref{p2} assume that that the triangulation of
 $W_{i+1}$ is chosen to be so small that the map $\gamma_i \circ
 \mu_{i+1}|_{...} : W_{i+1} \lo W_i$ admits a simplicial approximation
 $\omega_{i+1} : W_{i+1} \lo W_i$.  

  \medskip
{\bf(8)}\qua  Now we will verify a technical property which will be used
  later.  Namely we will show that
  $\kappa_{i+1},\psi_{i+1}|_{W_{i+1}}, \omega_{i+1}$ and $
  \psi_{i}|_{W_{i}}$ combinatorially commute. This property is
  equivalent to
$$\omega_{i+1} ((\kappa_{i+1}\circ\psi_{i+1})^{-1}(\Delta) \cap
 W_{i+1} ) \subset \psi_i^{-1}(\Delta)\text{ for every simplex }\Delta\text
 { of }K_i.$$
  Let $\Delta$ be a simplex of $K_i$.  Note that
 $\psi_i^{-1}(\Delta) $ is a subcomplex of $M_i$ with respect to both
 the CW and the simplicial structures of $M_i$.  Then since
 $\mu_{i+1}$ is a simplicial approximation of $
 \nu_{i+1}\circ\psi_{i+1}$ and $\psi_i^{-1}(\Delta)$ is a simplicial
 subcomplex of $M_ i$ we have
\begin{align*} 
  \mu_{i+1}
  ((\nu_{i+1}\circ\psi_{i+1})^{-1}(\psi_i^{-1}(\Delta)) )&\subset
  (\nu_{i+1}\circ\psi_{i+1})((\nu_{i+1}\circ\psi_{i+1})^{-1}(\psi_i^{-1}(\Delta))
  )\\& =\psi_i^{-1}(\Delta)\end{align*}
 and therefore
 $$ \mu_{i+1}
       ((\nu_{i+1}\circ\psi_{i+1})^{-1}(\psi_i^{-1}(\Delta)) \cap
       W_{i+1} ) \subset \psi_i^{-1}(\Delta)\cap S_i.$$  
$$( \gamma_i \circ \mu_{i+1} )
 ((\nu_{i+1}\circ\psi_{i+1})^{-1}(\psi_i^{-1}(\Delta)) \cap W_{i+1} )
 \subset\gamma_i( \psi_i^{-1}(\Delta)\cap S_i).\leqno{\rm thus}$$
 By (4) $
 \gamma_i( \psi_i^{-1}(\Delta)\cap S_i) \subset
 \psi_i^{-1}(\Delta)\cap W_i$ and since $\psi_i^{-1}(\Delta)\cap W_i $
 is a subcomplex of $W_i$ and $\omega_{i+1}$ is a simplicial
 approximation of $\gamma_i \circ \mu_{i+1}|_{...} : W_{i+1} \lo W_i$
 we get
$$\omega_{i+1} ((\nu_{i+1}\circ\psi_{i+1})^{-1}(\psi_i^{-1}(\Delta))
 \cap W_{i+1} ) \subset\psi_i^{-1}(\Delta)\cap W_i \subset
 \psi_i^{-1}(\Delta).$$
 Recall that $\kappa_{i+1}= \psi_i \circ
 \nu_{i+1}$ and finally get that
$$\omega_{i+1} ((\kappa_{i+1}\circ\psi_{i+1})^{-1}(\Delta) \cap
 W_{i+1} ) \subset \psi_i^{-1}(\Delta)\text{ for every simplex
 }\Delta\text{ of }K_i.$$

{\bf(9)}\qua $L_{i+1}$ is constructed as follows.  Let $C$ be an $(n+1)$-cell
 of $M_{i+1}$.  By (ii) of Proposition \ref{p1} there is an
 $(n+1)$-simplex $\Delta$ of $K_{i+1}$ such that
 $C=\psi_{i+1}^{-1}(\Delta)$.  By (iv) of Proposition \ref{p1} there
 is an $n$-simplex $\Delta_C$ of $W_{i+1}$ contained in $\partial C$
 which is sent by $\psi_{i+1}$ to the interior of $\Delta$.  Thus for
 every $(n+1)$-cell $C$ of $M_{i+1}$ we choose an $n$-simplex
 $\Delta_C$ of $W_{i+1}$ which consisits of free boundary points of
 $C$.

  We say that an $n$-simplex $\Delta_{i+1}$ of $W_{i+1}$ is marked
  (for replacement) if $\Delta_{i+1}$ is not among the chosen
  simplexes $\Delta_C$ and $\Delta_{i+1}$ is mapped by $\omega_{i+1}$
  onto an $n$-simplex $\Delta_i$ of $W_i$ such that $\Delta_i$ was
  replaced by an $n$-dimensional ${\cal L}$-ball
  $B_{\Delta_i}=\beta_i^{-1}( \Delta_i)$ while constructing $L_i$, see
  (5).

Replace each marked simplex $\Delta_{i+1}$ by an $n$-dimensional
${\cal L}$-ball $B_{\Delta_{i+1}}$ which is a copy of $B_{\Delta_i}$
attached by the simplicial map of the boundaries of $\Delta_{i+1}$ and
$B_{\Delta_{i+1}}$ induced by $\omega_{i+1}$ and $\beta_i$.  Such a
replacement induces the natural map $\lambda_{\Delta_{i+1}} :
B_{\Delta_{i+1}} \lo B_{\Delta_i}= \beta_i^{-1}( \Delta_i) \subset
L_i$.

 Let $C$ be an $ (n+1)$-cell of $M_{i+1}$ and let $T_C $ be the union
of all the simplexes of $W_{i+1}$ which are contained in $C$ and
different from $\Delta_C$.  Then after the replacement of all marked
simplexes of $W_{i+1}$ we will obtain from $T_C$ an $n$-dimensional
${\cal L}$-ball $B_C$ and the map $\lambda_C : B_C \lo L_i$ induced by
$\omega_{i+1}$ and $\beta_i$ for the simplexes of $T_C$ which were not
replaced and by the maps $ \lambda_{\Delta_{i+1}} $ for the replaced
simplexes.  Now replace $\Delta_C$ by an ${\cal L}$-ball
$B_{\Delta_C}$ which is a copy of $B_C$ attached by the identity map
of the boundaries and extend $ \lambda_C$ over $B_{\Delta_C}$ first
sending the points of $B_{\Delta_C}$ to the corresponding points of
$B_C$ and from there to $L_i$ using the map $ \lambda_C $. Thus we
construct $L_{i+1}$ and the map $\lambda_{i+1} : L_{i+1} \lo L_i$
having the property that for every $(n+1)$-cell $C$ of $M_{i+1}$ the
map
$$\lambda_{i+1} |_{...} : \beta^{-1}_{i+1}(\partial C) \lo
\lambda_{i+1}(\beta^{-1}_{i+1}(\partial C) )$$
 factors through an
$n$-dimensional ${\cal L}$-ball (namely through $B_{C})$.  

\medskip
{\bf(10)}\qua Define $\alpha_i = \psi_i \circ \beta_i : L_i \lo K_i$.  
  
Then the last property of (9) says that for every $(n+1)$-simplex
  $\Delta$ of $K_{i+1}$ the map $\lambda_{i+1} |_{...}  :
  \alpha_{i+1}^{-1}(\Delta) \lo \lambda_{i+1}
  (\alpha_{i+1}^{-1}(\Delta) )$ factors through a $\z_{({\cal L})}$-
  acyclic space.

\medskip
{\bf(11)}\qua It follows from (9) that if $a \in \Delta_C$ then $(\beta_i
   \circ \lambda_{i+1}) (\beta_{i+1}^{-1}(a)) \subset \omega_{i+1}(
   \partial C) $ and if $a\in W_{i+1}$ does not belong to any of the
   simplexes $ \Delta_C$ then simply $\beta_i \circ \lambda_{i+1}
   (\beta_{i+1}^{-1}(a)) = \omega_{i+1}( a)$.  

 Thus we obtain that
   for every simplex $\Delta$ of $K_{i+1}$
   $$(\beta_i \circ \lambda_{i+1})
   (\beta_{i+1}^{-1}(\psi_{i+1}^{-1}(\Delta)) = (\beta_i \circ
   \lambda_{i+1})(\alpha_{i+1}^{-1}(\Delta)) \subset \omega_{i+1}(
   \psi_{i+1}^{-1}(\Delta)\cap W_{i+1}) .$$

\medskip
{\bf(12)}\qua  Now we will show that the maps $\kappa_{i+1}, \alpha_{i+1},
\lambda_{i+1}$ and $ \alpha_i$ combinatorially commute.  Let $\Delta$
be a simplex of $K_i$.  Then $\kappa_{i+1}^{-1}(\Delta)$ is a
subcomplex of $K_i$ and by (11) and (8) we get that
\begin{align*}
 (\beta_{i}\circ\lambda_{i+1} )(( \kappa_{i+1}\circ
  \alpha_{i+1})^{-1}(\Delta))& = (\beta_{i}\circ\lambda_{i+1} )(
  \alpha_{i+1}^{-1}(\kappa_{i+1}^{-1}(\Delta))\\
  &\subset\omega_{i+1}( \psi_{i+1}^{-1}(\kappa_{i+1}^{-1}(\Delta))\cap
 W_{i+1}) \subset \psi_i^{-1}(\Delta).\end{align*}
 Thus we get that
$$(\alpha_{i}\circ\lambda_{i+1} )(( \kappa_{i+1}\circ
  \alpha_{i+1})^{-1}(\Delta)) \subset \Delta\text{ for every simplex }
  \Delta\text{ of }K_i$$
 and hence the maps $\kappa_{i+1}, \alpha_{i+1},
  \lambda_{i+1}$ and $ \alpha_i$ combinatorially commute.  

 Then by
  (iii) of Proposition \ref{p2} for every $a \in K_{i+1}$ and $b\in
  L_{i+1}$
 $$\lambda_{i+1}( \alpha_{i+1}^{-1}(st( a)) ) \subset
   \alpha_{i}^{-1}(st( \kappa_{i+1} (a)) )$$
$$\kappa_{i+1}(st((\alpha_{i+1} (b))) \subset st((\alpha_i\circ
   \lambda_{i+1})(b)).\leqno{\rm and}$$

{\bf(13)}\qua Let $a \in K_{i+1}$ be such that $\kappa_{i+1}(a)$ is not
        contained in the interior of any of the $(n+1)$-simplexes of
        $K_i$.  Then $\nu_{i+1}(a) \in W_i$.  Let us show that the map
        $\lambda_{i+1}|_{...} : \alpha_{i+1}^{-1}(st(a))) \lo
        \alpha_i^{-1}(st(\kappa_{i+1}(a)) )$ factors through a
        $\z_{({\cal L})}$- acyclic space.  

  By (4) there is a
        contractible subcomplex $N$ of $W_i$ such that
$$\gamma_i(st(\nu_{i+1}(a)) \cap S_i ) \subset N \subset
\psi_i^{-1}(st(\kappa_{i+1}(a))$$
 where $st(\nu_{i+1}(a))$ is the
star of $ \nu_{i+1}(a)$ with respect to the triangulation of $M_i$.

Since $\nu_{i+1}$ is combinatorial with respect to the
triangulation of $M_i$ we have
$$\nu_{i+1}(st(a)) \subset st(\nu_{i+1}(a))$$ and therefore
$$\nu_{i+1}(\psi_{i+1}(\psi_{i+1}^{-1}(st(a))) )= \nu_{i+1}(st(a))
 \subset st(\nu_{i+1}(a)).$$
 Since $\mu_{i+1}$ is a simplicial
 approximation of $\nu_{i+1} \circ \psi_{i+1}$ and $st(\nu_{i+1}(a))$
 is a simplicial subcomplex of $M_i$ we have
$$\mu_{i+1}(\psi_{i+1}^{-1}(st(a))) \subset st(\nu_{i+1}(a))$$ and
    therefore
$$\gamma_i(\mu_{i+1}(\psi_{i+1}^{-1}(st(a)) \cap S_i) \subset
       \gamma_i( st(\nu_{i+1}(a))\cap S_i) \subset N.$$  
 Then since
       $N$ is a subcomplex of $W_i$ and $\omega_{i+1}$ is a simplicial
       approximation of $\gamma_i \circ \mu_{i+1}|_{...} : W_{i+1} \lo
       W_i$ we have
$$ \omega_{i+1}(\psi_{i+1}^{-1}(st(a))\cap W_i) \subset N$$ and
hence
$$\omega_{i+1}(\psi_{i+1}^{-1}(st(a))\cap W_i) \subset N\subset
 \psi_i^{-1}(st(\kappa_{i+1}(a)).$$
 By (11)
$$\beta_i(\lambda_{i+1}(\alpha_{i+1}^{-1}(st(a)))) \subset
   \omega_{i+1}(\psi_{i+1}^{-1}(st(a))\cap W_i) $$ and therefore
$$\lambda_{i+1}(\alpha_{i+1}^{-1}(st(a)))\subset \beta_i^{-1}(N)
    \subset \beta^{-1}( \psi_i^{-1}(st(\kappa_{i+1}(a)))=
    \alpha_i^{-1}(st(\kappa_{i+1}(a))).$$
 By (5) the preimage under
    $\beta_i$ of every simplex of $W_i$ is $\z_{({\cal L})}$- acyclic
    and therefore by the combinatorial Vietoris-Begle theorem $
    \beta_i^{-1}(N)$ is $\z_{({\cal L})}$- acyclic.

  Thus the map $\lambda_{i+1}|_{...} : \alpha_{i+1}^{-1}(st(a))) \lo
\alpha_i^{-1}(st(\kappa_{i+1}(a)) )$ factors through a $\z_{({\cal
L})}$- acyclic space.

 \medskip
{\bf(14)}\qua  Let $a \in K_{i+2}$.  By (12) consider the following maps

{\leftskip 25pt $\lambda_{i+2}|_{...} : \alpha_{i+2}^{-1}(st( a)) \subset
   \alpha_{i+1}^{-1}(st( \kappa_{i+2} (a)) )$ and

  $\lambda_{i+1}|_{...} : \alpha_{i+1}^{-1}(st(\kappa_{i+2}( a))
   \subset \alpha_{i}^{-1}(st( \kappa_{i+1} (\kappa_{i+2}(a)) )$.\par}

 If $\kappa_{i+2}(a)$ is in the interior of an $(n+1)$-simplex
 $\Delta$ of $K_{i+1}$ then $st(\kappa_{i+2}( a))=\Delta$ and by (10)
 $\lambda_{i+2}|_{...} $ factors through a $\z_{({\cal L})}$- acyclic
 space.

 If $\kappa_{i+2}(a)$ is not in the interior of any of the
 $(n+1)$-simplexes of $K_{i+1}$ then by (13) $\lambda_{i+1}|_{...}$
 factors through a $\z_{({\cal L})}$- acyclic space.

 Thus we get that the map
$$\lambda_{i+1}\circ \lambda_{i+2}|_{...} : \alpha_{i+2}^{-1}(st( a))
   \lo \alpha_{i}^{-1}(st( \kappa_{i+1} (\kappa_{i+2}(a)) )$$
 always
   factors through a $\z_{({\cal L})}$- acyclic space.

   \medskip
{\bf(15)}\qua  Define
   
   {\leftskip 25pt
   $L=\invlim( L_i, \lambda_i)$ and $K=\invlim (K_i, \kappa_i)$ with
   the projections

   $\pi^L_i : L \lo L_i $ and $\pi^K_i : K \lo K_i $.\par}  

 It is clear
that $\dim L \leq n$.  For every $b=(b_i) \in L, b_i \in L_i$ define
the set $\pi(b)=\cap \{ (\pi_i^K)^{-1}(st(\alpha_i(b_i))):
i=1,2,3...\}\subset K$.  By (12) $\pi(b)\neq\emptyset$.  Recall that
by (6) diam($\kappa_i^j(\Delta)) <1/i$ for every simplex $ \Delta $ of
$K_i$ and $j <i$.  Then $\pi(b)$ is a singleton and by (12) one can
check that the function

{\leftskip 25pt
    $\pi : L \lo K$ defined by $b \lo \pi(b)$ is continuous with

    $\pi^{-1}(a) = \invlim ( \alpha_{i}^{-1}(st( a_i)), \lambda_i
   |_{...})$ for every $a=(a_i) \in K$.\par}  

 It is easy to see that for
   every simplex $\Delta$ of $K_i$, $\alpha_i^{-1}(\Delta) \neq
   \emptyset$.  Then we have that for every $a \in K$, $\pi^{-1}(a)$
   is not empty and by (14) we get that $\pi^{-1}(a)$ is $\z_{({\cal
   L})}$-acyclic.

 Thus $\pi$ is a $\z_{({\cal L})}$-acyclic surjective map.  

  \medskip
{\bf(16)}\qua 
 Let $ x \in X$ and $j \leq i$.  Define $f_{i}^j=
  \kappa_{i}^j \circ f_i : X \lo K_j$.  

 Then $f_{i+1}^i =
  \kappa_{i+1} \circ f_{i+1} = \psi_i \circ \nu_{i+1}\circ f_{i+1}$
  and by (6) $d_i(f_{i+1} (x), ( \psi_i\circ \phi_i)(x)) <
  \epsilon_i/2$.  

 By (2) and (i) of Proposition \ref{p1} we have
$$d_i(f_i(x) , ( \psi_i\circ \phi_i)(x))) < \epsilon_i/2\text{ and hence }
d_i(f_i(x), f_{i+1}^i (x)) < \epsilon_i.$$
 Then, since
 $f_{i+1}^j=\kappa_i^j \circ f_{i+1}$, (2) implies that $d_j(
 f_i^j(x), f_{i+1}^j (x)) < 1/4^i$.

Now define $h_j : X \lo K_j$ as $h_j =\lim_{i\to \infty} f_i^j$.

Then $d_j (h_j(x), f_j(x))< 1/2$.  Let $a=h_j(x)$ and let $y \in X$ be
such that $h_j(y)=h_j(x)$.  Then $d_j(f_j(y), f_j(x)) < 1$ and hence
by (1) $f_j(y) \in st^2( a) $ where $st^2(a)$ is the union of all
simplexes of $K_j$ intersecting $st(a)$.

Thus $h_j^{-1}(h_j(x)) \subset f_j ^{-1}(st^2(a))$ and hence by (6)
$h_j^{-1}(h_j(x))$ is of diam $\leq 4/j$.

It is clear that $h_{j}=\kappa_{j+1} \circ h_{j+1}$.  Then we have
that $h : X \lo K$ defined by $h(x)=(h_j(x)) \in K$ is an embedding of
$X$ in $K$.

   Identify $X$ with $h(X)$.  Then $r=\pi|_{...} : Z=\pi^{-1}(X) \lo
   X$ is the required resolution.  \endproof

 {\bf Remark}\qua
   If $n\geq 3$ then in (9) $B_{C}$ is simply connected and one can
   use the Hurewicz theorem to replace each simplex $\Delta_C$ by an
   $n$-cell attached to the boundary by a map of degree not divisible
   by the elements of $\cal L$ such that there is a map from the
   replacing cell to $B_{C}$.  In our construction we used an exact
   copy of $B_{C}$ for replacing $\Delta_C$ in order to guarantee that
   such a map exists for the non-simply connected case $n=2$.

  \end{section}\np

\Addresses\recd
\end{document}